\documentclass[a4paper]{amsart}
\usepackage[latin1]{inputenc}
\usepackage{amsmath}
\usepackage{amsthm}
\usepackage{amssymb}
\usepackage{amsfonts}
\usepackage{url}
\usepackage{apalike}

\newtheorem{theorem}{Theorem}[section]

\newtheorem{conj}[theorem]{Conjecture}

\newcommand\depth{\operatorname{depth}}
\newcommand\rank{\operatorname{rank}}
\newcommand\ideal[1]{\langle #1 \rangle}
\newcommand\iso\cong

\title{On low degree regular sequences in group cohomology}
\author{Mikael Johansson}
\address{
Mikael Johansson\\
Mathematisches Institut\\
Fakultät für Mathematik und Informatik\\
FSU Jena\\07743 Jena, Germany\\}
\email{mik@math.uni-jena.de}

\begin{document}

\begin{abstract}
  We investigate small $p$-groups with cohomology rings of depth
  higher than predicted by Duflot's theorem. In these groups, a
  sampling would suggest several naïve conjectures about the degrees of
  the additional regular sequence elements. We arrive at
  counterexamples to the idea that all cohomology rings exceeding
  Duflot's bound have degree 2 regular sequence elements past Duflot's
  bound as well as to the idea that regular sequences can be
  pulled back, preserving degrees, along a field extension.
\end{abstract}

\maketitle

\section{Introduction}

Let $G$ be a finite $p$-group. The cohomology ring $H^*(G,\mathbb
F_p)$ is a graded commutative, noetherian, local $\mathbb
F_p$-algebra. The depth of a graded commutative, noetherian, local
ring $(R,\mathfrak m,k)$ is defined as the largest $n$ such that there
are elements $x_1,\dots,x_n\in\mathfrak m$ such that $x_1$ is not a
zerodivisor in $R$, and $x_i$ is not a zerodivisor in
$R/\ideal{x_1,\dots,x_{i-1}}$ for all $i$. A sequence of $x_i$ in this
manner is called a \emph{regular sequence} in $R$, and we say that it
\emph{displays} the depth if the length is precisely the depth. We say
that the regular sequence $x_1,\dots,x_n$ has the \emph{degree
  sequence} $a_1,\dots,a_d$ if $|x_i|=a_i$ for any $i$.

Regular sequences displaying the depth of a ring $R$ are far from
unique. For the case $H^*(G,\mathbb F_p)$, a theorem by Duflot in
\cite{duflot81} states that $\depth H^*(G,\mathbb F_p)\geq \rank_p
Z(G)$. The proof of this theorem in \cite{duflot81} as well as the
proof in \cite{evens91}, both give an explicit regular sequence of
length $\rank_p Z(G)$.

We say that a group is \emph{superduflot} if its depth exceeds the
bound in Duflot's theorem. The Duflot regular sequence can be extended
to a maximal regular sequence, and we call the extending sequence a
superduflot regular sequence. An element with index higher than the
Duflot bound will be called a superduflot element.

For the groups referenced, we give the number in the small groups
library \cite{smallgroups}, and where applicable, the Hall-Senior
number \cite{hallsenior}

\section{Conjectures and counterexamples}

Initially, a small sampling of the calculated groups suggested the
conjecture \ref{conj:superduflot2}. While checking the conjecture
against further groups, the author was repeatedly able to find lower
superduflot degree sequences than given in the cohomology ring database
\cite{green_cohodb}. 

\begin{conj}\label{conj:superduflot2}
All cohomology rings of superduflot finite groups have superduflot
regular sequences in degrees at most 2
\end{conj}

Among the groups of order smaller than 256, the
calculations of \cite{carlson_cohodb} and \cite{green_cohodb} carry 
information about the cohomology rings of all $2$--groups up to order
64, and all $3$--groups and $5$--groups of order smaller than 243 as
well as 32 groups of order 128. In
total, there are 124 superduflot groups among these. 120 of these
admit regular sequences with the superduflot elements in degree
2.

Among the remaining 4 groups, an investigation of the to
date calculated rings revealed two groups of order 128 (small group
indices 929 and 934) and one group of order 243 (small groups index
28), all of which had a common annihilator for all cohomology classes
of degree lower than 4. Thus the conjecture \ref{conj:superduflot2}
cannot possibly hold.

Supposing, however, that the elements of low degree do not all have a
common annihilator. Then, we could note that $H^*(G,\mathbb
F_{p^r})\iso H^*(G,\mathbb F_p)\otimes_{\mathbb F_p}\mathbb F_{p^r}$
and use the linear independence of basis elements in $\mathbb F_{p^r}$
to construct an element over an extension field with the expected
degree and regularity properties. Thus, the statement refines into:

\begin{conj}
If $H^*(G,\mathbb F_{p^k})$ admits a regular sequence of length $d$
with degree sequence $a_1,\dots,a_d$, then $H^*(G,\mathbb F_{p})$
admits a regular sequence of equal length and with the same degree
sequence.
\end{conj}

On being shown this conjecture, Jon F. Carlson predicted that it would
be false for many $p$-groups. Unfortunately he was right, at least as
far as the validity of the conjecture is concerned. However there is
basically only one counterexample in the sample under investigation:

One group of order 64, namely small group index 139 or Hall-Senior
260, has depth 2. There is one regular element $q$ of degree 8 of the
cohomology ring, and it can be augmented with an element of degree 4
to a regular sequence displaying the depth. \cite{carlson_cohodb}
gives the regular sequence, using the presentation below, as

$$
q,z^4  + z^2 w + zx^3  + zxv + x^4  + x^2 v + w^2  + wv + v^2 
$$

We let the ring $H^*(G,\mathbb F_2)$ be presented as in
\cite{carlson_cohodb} with the generators $z,y,x$ in degree 1,
$w,v$ in degree 2, $u$ in degree 3, $t,s,r$ in degree 5 and $q$ in
degree 8, and with the relations generated by
\begin{align*}
&zy,yx,y^3,zv + xw,
zu + y^2 w, 	z^2 v + zxv, 	y^2 v + xu, 	zwv + zv^2, \\
&y^2 u, 	zxv^2 + zr + xt + xs + xr,
ywu + yvu + yr + w^2 v + wv^2 + u^2, \\
&zvu + xs + xr, 	ys + xvu, 	yt,
zs + zr + xs + xr, \\
&zv^3 + ywv^2 + yu^2 + ws + vt, 	y^2 r,
ywvu + ywr + xvs + xvr + ut, \\
&ywr + yv^2 u + yvs + yvr + ut + us, 	zwr + zvr,
yvu^2 + yur + wvt + v^2 t, \\
&z^5 r + z^4 w^3 + z^3 wt + z^2 q + zx^2 vr + zv^2 r + x^3 vt + xv^2 t
+ xv^2 s + xv^2 r + t^2, \\
&z^3 x^2 r + zv^2 r + yv^2 r + x^3 vt + x^2 v^4 + x^2 q + xv^2 t + vut + sr,\\
&z^4 xr + zxq + zv^2 r + x^3 vt + xv^2 s + xv^2 r + ts, \\
&y^2 q + yv^2 r + xv^2 s + xv^2 r + wvu^2 + ts + tr + sr + r^2, \\
&ywvr + wus + vut + ts + tr,
yv^2 s + yv^2 r + xv^2 s + xv^2 r + vut + s^2 + sr 
\end{align*}

A calculation with Singular \cite{singular} now shows that
$\xi_1=w+x^2$ and $\xi_2=v+y^2$ have disjoint annihilators, and
furthermore that with $\omega$ the additional basis element, that
$\xi_1+\omega\xi_2$ is regular in $H^*(G,\mathbb F_4)/\ideal{q}$. Thus
there is a regular sequence in $H^*(G,\mathbb F_{2^2})$ with degree
sequence $8,2$.

Now, the annihilator of $y^2$ in $H^*(G,\mathbb F_2)$ is generated by
the classes $x,y,z,u,r,s,t$ and $wv^2+w^2v$. Thus $y^2wv, y^2w(w+v)$
and $y^2v(w+v)$ all are non-zero, and annihilate $w+v,v$ and $w$
respectively. Thus, every element of $H^*(G,\mathbb F_2)$ with degree
less than 4 is annihilated by one of these three classes.

\bibliographystyle{mik}
\bibliography{library}

\begin{thebibliography}{7}
\providecommand{\natexlab}[1]{#1}
\providecommand{\url}[1]{\texttt{#1}}
\providecommand{\urlprefix}{URL }
\expandafter\ifx\csname urlstyle\endcsname\relax
  \providecommand{\doi}[1]{doi:\discretionary{}{}{}#1}\else
  \providecommand{\doi}{doi:\discretionary{}{}{}\begingroup
  \urlstyle{rm}\Url}\fi
\providecommand{\bibAnnoteFile}[1]{%
  \IfFileExists{#1}{\begin{quotation}\noindent\textsc{Key:} #1\\
  \textsc{Annotation:}\ \input{#1}\end{quotation}}{}}
\providecommand{\bibAnnote}[2]{%
  \begin{quotation}\noindent\textsc{Key:} #1\\
  \textsc{Annotation:}\ #2\end{quotation}}
\providecommand{\eprint}[2][]{\url{#2}}

\bibitem[Besche et~al., 2001]{smallgroups}
Hans~Ulrich Besche, Bettina Eick, and E.~A. O'Brien.
\newblock The groups of order at most 2000.
\newblock \emph{Electron. Res. Announc. Amer. Math. Soc.}, 7:1--4 (electronic)
  2001.
\newblock ISSN 1079-6762.
\bibAnnoteFile{smallgroups}

\bibitem[Carlson, 2001]{carlson_cohodb}
Jon~F. Carlson.
\newblock Cohomology computations 2001.
\newblock \urlprefix\url{http://www.math.uga.edu/~lvalero/cohointro.html}.
\bibAnnoteFile{carlson_cohodb}

\bibitem[Duflot, 1981]{duflot81}
J.~Duflot.
\newblock Depth and equivariant cohomology.
\newblock \emph{Comment. Math. Helv.}, 56(4):627--637 1981.
\newblock ISSN 0010-2571.
\bibAnnoteFile{duflot81}

\bibitem[Evens, 1991]{evens91}
Leonard Evens.
\newblock \emph{The cohomology of groups}.
\newblock Oxford Mathematical Monographs. Oxford University Press, New York
  1991.
\newblock ISBN 0-19-853580-5.
\bibAnnoteFile{evens91}

\bibitem[Green, 2006]{green_cohodb}
David~J. Green.
\newblock Cohomology rings of finite groups 2006.
\newblock
  \urlprefix\url{http://www.math.uni-wuppertal.de/~green/Coho_v2/index.html}.
\bibAnnoteFile{green_cohodb}

\bibitem[Greuel et~al., 2001]{singular}
G.-M. Greuel, G.~Pfister, and H.~Schönemann.
\newblock Singular, 2-0-3. a computer algebra system for polynomial
  computations. 2001.
\newblock \urlprefix\url{http://www.singular.uni-kl.de}.
\bibAnnoteFile{singular}

\bibitem[Hall and Senior, 1964]{hallsenior}
Marshall Hall, Jr. and James~K. Senior.
\newblock \emph{The groups of order {$2\sp{n}\,(n\leq 6)$}}.
\newblock The Macmillan Co., New York 1964.
\bibAnnoteFile{hallsenior}

\end{thebibliography}


\end{document}